\documentclass[11pt]{amsart}
\usepackage{amsmath,amsfonts,amsthm,amssymb}

\newtheorem{thm}{Theorem}
\newtheorem{lem}[thm]{Lemma}

\newtheorem{corol}[thm]{Corollary}

\renewcommand{\i}{\infty}
\renewcommand{\a}{\alpha}

\renewcommand{\l}{\lambda}
\renewcommand{\t}{\tau}

\newcommand{\rt}{\right}
\newcommand{\lt}{\left}
\numberwithin{equation}{section}
\numberwithin{thm}{section}
\begin{document}

\title{Partitions weighted by the parity of the crank}
\author{Dohoon Choi, Soon-Yi Kang, and Jeremy Lovejoy}

\date{\today}

\address{Korea Institute for Advanced Study, 207-43 Cheongnyangni-2-dong, Dongdaemun-gu, Seoul, 130-722, KOREA}
\email{choija@kias.re.kr}
\address{Korea Institute for Advanced Study, 207-43 Cheongnyangni-2-dong, Dongdaemun-gu, Seoul, 130-722, KOREA}
\email{sykang@kias.re.kr}
\address{LIAFA, CNRS and Universit\'e Denis Diderot,
2, Place Jussieu, Case 7014, F-75251 Paris Cedex 05, FRANCE}
\email{lovejoy@liafa.jussieu.fr}

\thanks{The third author was partially supported by an ACI
``Jeunes Chercheurs et Jeunes Chercheuses"}

\begin{abstract}  A partition statistic \lq crank' gives combinatorial
interpretations for Ramanujan's famous partition congruences. In
this paper, we establish an asymptotic formula, Ramanujan type
congruences, and q-series identities that the number of partitions
with even crank $M_e(n)$ minus the number of partitions with odd
crank $M_o(n)$ satisfies. For example, we show that
$M_e(5n+4)-M_o(5n+4)\equiv 0 \pmod 5.$ We also determine the exact
values of $M_e(n)-M_o(n)$ in case of partitions into distinct
parts, which are at most two and zero for infinitely many $n$.

\end{abstract}
\maketitle

\section{Introduction}\label{intro}
A partition $\l$ of a positive integer $n$ is a weakly decreasing
sequence of positive integers $\l_1\geq\l_2\geq\cdots\geq\l_k$
such that $\l_1+\l_2+\cdots+\l_k=n$. Then $\l_1$ is the largest
part and $k$ is the number of parts of $\l$ .  In 1944, Dyson
\cite{dyson} defined the \textit{rank of a partition $\l$} by
$$\mathrm{rank}(\l):=\l_1-k$$
and observed that this partition statistic appeared to give
combinatorial interpretations for Ramanujan's famous congruences
$p(5n+4)\equiv 0\pmod 5$ and $p(7n+5)\equiv 0\pmod 7$, where
$p(n)$ is the number of partitions of $n$.  Dyson also conjectured
the existence of another partition statistic, named the \lq
crank', that would explain the last Ramanujan's partition
congruence $p(11n+6)\equiv 0\pmod {11}$ combinatorially.  His
observations on the rank were first proved by Atkin and
Swinnerton-Dyer \cite{ASD} in 1954 and the crank was found by
Andrews and Garvan \cite{AG} in 1988.  This Andrews-Garvan crank
is defined by
$$\mathrm{crank}(\l)=\left\{%
\begin{array}{ll}
    \l_1, & \hbox{{if} $\mu(\l)=0$;} \\
    \nu(\l)-\mu(\l), & \hbox{{if} $\mu(\l)>0$,} \\
\end{array}%
\right.$$ where $\mu(\l)$ denotes the number of ones in $\l$ and
$\nu(\l)$ denotes the number of parts of $\l$ that are strictly
larger than $\mu(\l)$.

Studying the number of partitions with even rank minus the number
with odd rank has led to some rather intriguing mathematics. For
example, if we let $N_e(n)$ (resp. $N_o(n)$) denote the number of
partitions of $n$ with even (resp. odd) rank, then we have
\begin{equation}
\sum_{n=0}^\i(N_e(n)-N_o(n))q^n=\sum_{n=0}^\i\frac{q^{n^2}}{(-q;q)_n^2}=:f(q),
\end{equation}
where
$$(z;q)_n=\prod_{j=0}^{n-1}(1-zq^j).$$
This $q$-series is one of Ramanujan's third order mock theta
functions and has been the subject of a number of works (e.g.
\cite{Andrewsmock,AL,BO,BO2,Dragonette,Lewis,Watson}).  Most
recently, Bringmann and Ono \cite{BO,BO2} have shown that $f(q)$
belongs to the theory of \emph{weak Maass forms}, and they used
this framework to prove a longstanding conjecture of Andrews and
Dragonette, giving an exact formula for $N_e(n) - N_o(n)$:

\begin{thm} \cite[Theorem 1.1]{BO}\label{form-rank} Let $I_s(x)$ be the usual I-Bessel function of order $s$.
 If $n$ is a positive integer, then
$$N_e(n)-N_o(n)=\pi(24n-1)^{-1/4}\sum_{k=1}^\i
\frac{(-1)^{\lfloor(k+1)/2\rfloor}A_{2k}(n-\frac{k(1+(-1)^k}{4})}{k}I_{1/2}\lt(\frac{\pi\sqrt{24n-1}}{12k}\rt),$$
where
$$A_k(n):=\frac{1}{2}\sqrt{\frac{k}{12}}
\mathop{\sum_{x (\mathrm{mod} 24k)}}_{x^2\equiv -24n+1
(\mathrm{mod} 24k)}\chi_{12}(x)e^{2\pi i x/12k},$$ with the sum
running over the residue classes modulo $24k$ and
$\displaystyle{\chi_{12}(x):=\lt(\frac{12}{x}\rt)}$.
\end{thm}

For another example, let $N_e(\mathcal{D},n)$ (resp.
$N_o(\mathcal{D},n)$) denote the number of partitions into
distinct parts with even (resp. odd) rank. Then
$$\sum_{n=0}^\i(N_e(\mathcal{D},n)-N_o(\mathcal{D},n))q^n=\sum_{n=0}^\i\frac{q^{n(n+1)/2}}{(-q;q)_n}=:R(q).$$
In \cite{ADH}, Andrews, Dyson, and Hickerson showed that the
coefficients of $R(q)$ have multiplicative properties determined
by a certain Hecke character associated to the ring of integers of
the real quadratic field $\mathbb{Q}(\sqrt 6)$, and Cohen
\cite{Cohen} subsequently showed that $R(q)$ belongs to the theory
of \emph{Maass waveforms}. To state the main theorem of
\cite{ADH}, we must note that every number which is $1$ modulo $6$
and greater than $1$ has a factorization of the form
\begin{equation}\label{factor} 6m+1=p_1^{e_1}p_2^{e_2}\cdots
p_r^{e_r},\end{equation} where $r\geq 1$, each $p_i$ is either a
prime $\equiv 1 \pmod 6$ or the negative of a prime $\equiv 5
\pmod 6$, the $p_i$'s are distinct, and $e_i$'s are positive
integers.

\begin{thm}\cite[Theorem 3]{ADH}\label{ADHthm}
For all positive integers $n$ we have
$N_e(\mathcal{D},n)-N_o(\mathcal{D},n)=T(24n+1)$, where $24n+1$
has the factorization (\ref{factor}) and $T(m)=T(p_1^{e_1})\cdots
T(p_r^{e_r})$ is the multiplicative function defined on powers of
primes by:
 $$ T(p^e) =
\begin{cases}
0, & \text{if $p \not\equiv 1 \pmod{24}$ and $e$ is odd},\\
1, & \text{if $p \equiv 13\ \mathrm{or}\ 19 \pmod{24}$ and $e$ is even}, \\
(-1)^{e/2}, & \text{if $p \equiv 7 \pmod{24}$ and  $e$ is even}, \\
e+1, & \text{if $p \equiv 1 \pmod{24}$ and  $T(p)=2$}, \\
(-1)^e(e+1), & \text{if $p \equiv 1 \pmod{24}$ and $T(p)=-2$}.
\end{cases}$$
\end{thm}

This ``exact" formula has numerous interesting consequences, such
as the fact that $N_e(\mathcal{D},n) - N_o(\mathcal{D},n)$ is
almost always zero and assumes every integer infinitely often.

In this paper we pass from the rank to the crank, studying the
number of partitions with even crank minus the number with odd
crank.  In doing so, we leave the world of weak Maass forms and
Maass waveforms and enter the world of classical modular forms.
Specifically, if $M(m,n)$ denotes the number of partitions of $n$
with crank $m$, then the generating function for $M(m,n)$ is given
in \cite{AG} by

\begin{equation}\label{gencrank}
\sum_{m=-\i}^\i\sum_{n=0}^\i
M(m,n)a^mq^n=\frac{(q;q)_\i}{(aq;q)_\i(q/a;q)_\i},\end{equation}
where
$$(z;q)_\i=\lim_{n \to \i}(z;q)_n.$$
Letting $M_e(n)$ (resp. $M_o(n)$) denote the number of partitions
of $n$ with even (resp. odd) crank, then by setting $a=-1$ in
(\ref{gencrank}) we have

\begin{equation}\label{g}\sum_{n=0}^\i(M_e(n)-M_o(n))q^n=\frac{(q;q)_\i}{(-q;q)_\i^2}:=g(q).\end{equation}

Apparently the only study of $g(q)$ was done by Andrews and Lewis
\cite{AL}, who proved that $M_e(n)>M_o(n)$ if $n$ is even and
$M_e(n)<M_o(n)$ if $n$ is odd (by showing all the coefficients of
$g(-q)$ are positive).  Here we shall examine several other
aspects of $g(q)$, such as congruences, asymptotics, and
$q$-series identities.

We begin with congruence properties satisfied by $M_e(n) -
M_o(n)$. From work of Treneer \cite{Tr1}, we know immediately that
$M_e(n) - M_o(n)$ has infinitely many congruences in arithmetic
progressions modulo any prime coprime to $6$.  The obvious
question, then, is whether any of these congruences are as simple
and elegant as those of Ramanujan for the partition function.  We
prove in fact that the crank difference $M_e(n)-M_o(n)$ satisfies
a family of congruences modulo powers of $5$.

\begin{thm} \label{family}
For all $\alpha\geq 0$, we have
$$ M_e(n) - M_o(n) \equiv 0 \pmod{5^{\alpha+1}}\quad \mathrm{if}\quad 24n\equiv 1
\pmod{5^{2\alpha+1}}.$$
\end{thm}

Following from the proof of Theorem \ref{family} will be a
generating function for $M_e(5n+4) - M_o(5n+4)$:
\begin{thm} \label{ramatype}
$$\sum_{n=0}^{\infty} (M_e(5n+4) - M_o(5n+4))q^n =
\frac{5(q;q^2)_{\infty}^2(q^5;q^5)_{\infty}(q^{10};q^{10})_{\infty}^2}{(q^2;q^2)_{\infty}^2}.
$$
\end{thm}

Next we apply the Hardy-Ramanujan circle method to the generating
function $g(q)$, obtaining an asymptotic formula for
$M_e(n)-M_o(n)$:

\begin{thm}\label{asy-cr} If $n$ is a positive integer, then
$$M_e(n)-M_o(n)=\frac{1}{\sqrt{n-1/24}}\sum_{0<k<5\sqrt{n}/2}
\frac{B_k(n)}{\sqrt{k}}\cosh\lt(\frac{\pi}{k}\sqrt{\frac{n-1/24}{6}}\rt)+E_n,$$
where
$$B_k(n)=\sum_{(h,2k)=1}e^{\pi i (2s(h,k)-3s(h,2k))}e^{-2\pi i n h/2k},$$
for the Dedekind sum $s(h,k)$, and where $|E_n|<194
n^{1/4}.$\end{thm}

To finish our study of the function $M_e(n) - M_o(n)$, we discuss
consequences of some $q$-series expansions for $g(q)$.  As an
example, we give a weighted partition identity involving $M_e(n) -
M_o(n)$.  To state this we use the notion of the ``initial run" of
a partition, by which we mean the largest increasing sequence of
part sizes starting with $1$. For example, the partition
$(7,7,5,3,3,3,3,2,1,1)$ has initial run $(1,2,3)$, while the
partition $(6,6,5,2,2,2,2)$ has no initial run at all.
\begin{thm} \label{weighted}
For a nonempty partition $\lambda$, define the weight
$\omega(\lambda)$  to be
$$
\omega(\lambda) = 1 + 4\sum_{j} (-1)^j,
$$
where the sum is over those $j$ in the initial run which occur an
odd number of times in $\lambda$.  Then
$$
M_e(n) - M_o(n) = \sum_{\lambda} \omega(\lambda),
$$
where the sum is over all partitions $\lambda$ of $n$.
\end{thm}

For example, take $n=4$. The partition $4$ has weight $1$, $(3,1)$
has weight $1-4 = 3$, $(2,2)$ has weight $1$, $(2,1,1)$ has weight
$1+4 = 5$, and $(1,1,1,1)$ has weight $1$. Summing these weights
gives $1 - 3 + 1 + 5 + 1 = 5$, which is, as expected, $M_e(4) -
M_o(4)$.

Finally, inspired by the work of Andrews, Dyson and Hickerson in
Theorem \ref{ADHthm}, we look at what happens if we restrict our
crank difference to partitions into distinct parts. In this case,
the definition of the crank simplifies considerably, and we are
able to use basic manipulations of $q$-series to prove an exact
formula.

Let $M_e(\mathcal{D},n)$ (resp. $M_o(\mathcal{D},n)$) denote the
number of partitions into distinct parts with even (resp. odd)
crank.  For partitions into distinct parts the crank is either the
largest part, if there is no one appearing, or the number of parts
minus $2$ if there is a one.  Let $\mathcal{P}$ denote the set of
pentagonal numbers, i.e., numbers of the form $m(3m+1)/2$ for $m$
an integer. If $n = m(3m+1)/2$, we write $R(n) = m$.  Finally, we
use the notations $\lfloor n \rfloor_p$ and $\lceil n \rceil_p$,
to denote the pentagonal floor and ceiling of $n$, i.e., the
largest (resp. smallest) pentagonal number $\leq$ (resp. $\geq$)
$n$.
\begin{thm} \label{formula}
For positive integers $n$ we have
$$ M_e(\mathcal{D},n) - M_o(\mathcal{D},n) =
\begin{cases}
1, & \text{if $n \in \mathcal{P}$ and $R(n)$ is odd and positive},
\\
-1, & \text{if $n \in \mathcal{P}$ and $R(n)$ is not as above}, \\
2, & \text{if $n \not \in \mathcal{P}$, $R(\lfloor n \rfloor_p)$
is odd and positive, and $n \equiv \lfloor n \rfloor_p \pmod{2}$}, \\
-2, & \text{if $n \not \in \mathcal{P}$, $R(\lfloor n \rfloor_p)$
is even and positive, and $n \equiv \lfloor n \rfloor_p
\pmod{2}$}, \\
-2(-1)^{n - \lfloor n \rfloor_p}, & \text{if $n \not \in
\mathcal{P}$ and $R(\lfloor n
\rfloor_p)$ is even and negative}, \\
0, & \text{otherwise}.
\end{cases}
$$
\end{thm}

\begin{corol}
The quantities $M_e(\mathcal{D},n)$ and $M_o(\mathcal{D},n)$
differ by at most $2$ and are equal for infinitely many $n$.
\end{corol}

To illustrate the above theorem, take $n=6$.  Then $n$ is not a
pentagonal number and $R(\lfloor 6 \rfloor_p) = R(5) = -2$.  Hence
by the penultimate case we expect $M_e(\mathcal{D},6) -
M_o(\mathcal{D},6)$ to be $-2(-1)^{6-5} = 2$.  Indeed, the four
partitions of $6$ into distinct parts are $(6)$, $(5,1)$, $(4,2)$,
and $(3,2,1)$, each of which has even crank except for $(3,2,1)$.

The paper is organized as follows.  In the next section, we prove
the family of congruences in Theorem \ref{family} and the
generating function in Theorem \ref{ramatype}.  In Section $3$, we
establish the asymptotic formula in Theorem \ref{asy-cr}.  In
Section $4$, we prove the weighted identity in Theorem
\ref{weighted} and discuss a similar identity.  Finally, In
Section $5$ we treat the exact formula in Theorem \ref{formula}.

\section{Proof of Theorems \ref{family} and \ref{ramatype}}
Here we follow the exposition in Gordon and Hughes' \cite{Go-Hu1}
rediscovery of some congruences of R\o dseth \cite{Rodseth}. The
reader may should have some familiarity with the preliminaries in
\cite{Go-Hu1}. Define
$$F(\t):=\frac{\eta(\t)^3}{\eta(2\t)^2}\frac{\eta(50\t)^2}{\eta(25\t)^3}.$$
Here $\eta(\tau) := q^{1/24} \prod_{n \geq 1} (1-q^n)$ and
$q:=e^{2 \pi i \t}$.  Applying \cite[Theorem 2, Theorem
3]{Go-Hu1}, we have $F(\t) \in M_0(\Gamma_0(50))$, the space of
weight $0$ modular functions on $\Gamma_0(50)$. Then $F(\t)$ is
holomorphic on the upper half plane $\mathbb{H}$ and its orders at
the cusps $\nu/\delta$ are as given below by \cite[Theorem
3]{Go-Hu1}.
\begin{center}
\begin{tabular}{|c|c|c|c|c|c|c|}
  \hline
  $\delta$ & 1 & 2 & 5 & 10 & 25 & 50 \\
  \hline
  $ord_{\nu/\delta}F$
 & 4 & -1 & 0 & 0 & -4 & 1 \\
  \hline
\end{tabular}
\end{center}
Next, recall the $U_d$-operator, which acts on power series by
$$
\sum_{n=0}^{\infty} a(n)q^n \big | U_d := \sum_{n=0}^{\infty}
a(dn)q^n.
$$ Note that if $f$ and $g$ are power series in $q$, we have
$$
\left(f(q^d)g(q)\right) \big | U_d = f(q) \left(g(q) \big | U_d
\right).
$$
By \cite[Theorem 5]{Go-Hu1} we have that $F(\t)|U_5\in
M_0(\Gamma_0(10))$ with the following lower bounds for the orders
of $F|U_5$ at the cusps.
\begin{center}\begin{tabular}{|c|c|c|c|c|}
  \hline
 $\delta$ & 1 & 2 &5 & 10 \\
  \hline
   $ord_{\nu/\delta}F|U_5\geq$ & 0 & -1 & 0 & 1 \\
  \hline
\end{tabular}
\end{center}

Now consider the function
$$G(\t):=\frac{\eta(\t)^2}{\eta(2\t)^4}\frac{\eta(10\t)^4}{\eta(5\t)^2}.$$
Applying \cite[Theorem 2, Theorem 3]{Go-Hu1}, we find that
$G(\t)\in M_0(\Gamma_0(10))$ and that its orders at the cusps are
as follows.
\begin{center}\begin{tabular}{|c|c|c|c|c|}
  \hline
 $\delta$ & 1 & 2 &5 & 10 \\
  \hline
   $ord_{\nu/\delta}G$ & 0 & -1 & 0 & 1 \\
  \hline
\end{tabular}
\end{center}
Since the only holomorphic modular functions of weight $0$ are the
constant functions, comparing the last two tables and the Fourier
series expansions of $F|U_5$ and $G$ gives $F|U_5=5G$.  Now, if we
consider $G(\t)$ as a function in $M_0(\Gamma_0(50))$, rather than
the subfield $M_0(\Gamma_0(10))$, we find from \cite[Theorem
3]{Go-Hu1} that its orders at the cusps are as follows.
\begin{center}
\begin{tabular}{|c|c|c|c|c|c|c|}
  \hline
  $\delta$ & 1 & 2 & 5 & 10 & 25 & 50 \\
  \hline
  $ord_{\nu/\delta}G$
 & 0 & -5 & 0 & 1 & 0 & 1 \\
  \hline
\end{tabular}
\end{center}
Hence by  \cite[Theorem 5]{Go-Hu1} applied to $G^i$ and $FG^i$,
these functions are on $\Gamma_0(10)$ and we have the following
lower bounds for the orders of $G^i|U_5$ and $FG^i|U_5$ at the
cusps.
\begin{center}
\begin{tabular}{|c|c|c|c|c|}
  \hline
  $\delta$ & 1 & 2 & 5 & 10 \\
  \hline
 $ord_{\nu/\delta}G^i|U_5\geq$& 0 & $-5i$ & 0 & $i/5$ \\
  $ord_{\nu/\delta}FG^i|U_5\geq$ &0 & $-5i-1$ & 0 & $(i+1)/5$ \\
  \hline
\end{tabular}\end{center}
If $i\geq 0$, this implies that $G^i|U_5$ and $FG^i|U_5$ are
polynomials in $G$ of degrees at most $5i$ and $5i+1$
respectively. Hence \begin{equation}\label{series}
G^i|U_5=\sum_{j\geq 0}a_{ij}G^j\quad \mathrm{and}\quad
FG^i|U_5=\sum_{j\geq 0}b_{ij}G^j\end{equation} for complex
coefficients $a_{ij}$ and $b_{ij}$.

Let $S$ be the vector space of all polynomials $P=\sum_{j\geq
0}c_jG^j$ and $T$ be the subspace of such polynomials with $0$ as
constant terms.  By considering our lower bounds for the orders of
$G^i|U_5$ and $FG^i|U_5$ and (\ref{series}), we see that $U_5$
maps $S$ to itself as well as $T$ to itself.  In addition, the
linear transformation $V:P\to (FP)|U_5$ maps $S$ into $T$. With
respect to the basis $G$, $G^2$, $G^3$, $\dots$ of $T$ the
matrices of $U_5$ and $V$ restricted to $T$ are respectively
$$U_5=(A:=(a_{ij}))\quad \mathrm{and }\quad  V=(B:=(b_{ij})),$$
for $1\leq i,j<\i.$

If we define a sequence of functions $L_\nu$ $(\nu\geq 0)$
inductively by putting for $\a\geq 0,$

\begin{equation}\label{defL} L_0=1,\ L_{2\a+1}=FL_{2\a}|U_5,\ \mathrm{and} \
L_{2\a+2}=L_{2\a+1}|U_5,
\end{equation}
then $$L_1=5G=(5,0,0,\dots),$$ $$L_{2\a+1}=(5,0,0,\dots)(AB)^\a,$$
and
$$L_{2\a+2}=(5,0,0,\dots)(AB)^\a A.$$

On the other hand, it follows from induction on $\a$ that
\begin{equation}\label{claimL}
L_{2\a+1}=\frac{(q^{10};q^{10})^2_\i}{(q^{5};q^{5})^3_\i}\sum_{n=1}^\i
(M_e(m)-M_o(m))q^n,\end{equation} where
$m=5^{2\a+1}n-1-5^2-\cdots-5^{2\a}$.
Theorem \ref{family} will then follow from the following theorem.

\begin{thm}\label{main}If we set
$L_{2\a+1}=(l_1(2\a+1),l_2(2\a+1), \dots),$ then $l_i(2\a+1)$ are
integers divisible by $5^{\a+1}$.\end{thm}

We prove Theorem \ref{main} in three steps.

\begin{lem}\label{a} For all $i$, $j$, we have $a_{ij}\in \mathbb{Z}$ and
$$\pi(a_{ij})\geq [\frac{5j-i-1}{6}],$$
where $\pi(n)$ denotes the 5-adic order of $n$.\end{lem}
\begin{proof}
Define
$$\phi(\t):=\frac{\eta(\t)}{\eta(2\t)^2}\frac{\eta(50\t)^2}{\eta(25\t)}\in
M_0(\Gamma_0(50)).$$ Again, by \cite[Theorem 3, Theorem
5]{Go-Hu1}, we obtain the orders of $\phi$ and the lower bounds
for the orders of $\phi^\mu|U_5$ for any nonnegative number $\mu$
at the cusps as follows.
\begin{center}
\begin{tabular}{|c|c|c|c|c|c|c|}
  \hline
  $\delta$ & 1 & 2 & 5 & 10 & 25 & 50 \\
  \hline
  $ord_{\nu/\delta}\phi$
 & 0 & -3 & 0 & 0 & 0 & 3 \\
  \hline
\end{tabular}
\end{center}

and

\begin{center}\begin{tabular}{|c|c|c|c|c|}
  \hline
 $\delta$ & 1 & 2 &5 & 10 \\
  \hline
   $ord_{\nu/\delta}\phi^\mu|U_5\geq$ & 0 & -3$\mu$ & 0 & 3$\mu/5$
    \\
  \hline
\end{tabular}
\end{center}

Hence for any $\mu\geq 0$,
$$\phi(\t)^\mu|U_5=\frac{1}{5}\sum_{\lambda=0}^4\phi(\frac{\t+\lambda}{5})^\mu\in
M_0(\Gamma_0(10))$$ is a polynomial in $G$ of degree at most
$3\mu$. Then by Newton's identities relating power sums of the
roots to the coefficients of a polynomial (see \cite{Mead} for
example), $\phi(\frac{\t+\lambda}{5})$ $(0\leq\lambda\leq 4)$ are
the roots of an equation
\begin{equation}\label{newton}t^5-\sigma_1 t^4+\sigma_2 t^3-\sigma_3 t^2+\sigma_4 t-\sigma_5=0,\end{equation}
where $\sigma_i's$ are elementary symmetric functions in
$\mathbb{C}[G]$.  As a matter of convenience, we consider the
reciprocal equation
$$u^5-\frac{\sigma_4}{\sigma_5}u^4
+\frac{\sigma_3}{\sigma_5}u^3-\frac{\sigma_2}{\sigma_5}u^2+\frac{\sigma_1}{\sigma_5}u-\frac{1}{\sigma_5}=0.$$
The roots of this equation are the functions
$\phi(\frac{\t+\lambda}{5})^{-1}$ $(0\leq\lambda\leq 4)$.

Note that
\begin{center}
\begin{tabular}{|c|c|c|c|c|c|c|}
  \hline
  $\delta$ & 1 & 2 & 5 & 10 & 25 & 50 \\
  \hline
  $ord_{\nu/\delta}\phi^{-1}$
 & 0 & 3 & 0 & 0 & 0 & $-3$ \\
  \hline
\end{tabular}
\end{center}

and

\begin{center}\begin{tabular}{|c|c|c|c|c|}
  \hline
 $\delta$ & 1 & 2 &5 & 10 \\
  \hline
   $ord_{\nu/\delta}\phi^{-\mu}|U_5\geq$ & 0 & 0 & 0 & $-3\mu/5$
    \\
  \hline
\end{tabular}
\end{center}

Since both $\phi^{-\mu}|U_5$ and $G^{-1}$ belong to $
M_0(\Gamma_0(10))$ and are holomorphic except for a pole at
infinity, it follows that $\phi^{-\mu}|U_5$ is a polynomial in
$G^{-1}$ of degree at most $[3u/5]$. From
$$G^{-1}(\t)=q^{-1}+2+q+2q^2+\cdots$$ and
$$G^{-2}(\t)=q^{-2}+4q^{-1}+6\cdots,$$
we can find that
$$\phi^{-1}|U_5=1,$$
$$\phi^{-2}|U_5=2q^{-1}+3+O(q)=2G^{-1}-1,$$
$$\phi^{-3}|U_5=6q^{-1}+7+O(q)=6G^{-1}-5,$$
$$\phi^{-4}|U_5=6q^{-2}+24q^{-1}+31+O(q)=6G^{-2}-5,$$
which happen to give exactly same polynomials in equations (17) in
\cite{Go-Hu1}.  Now the rest of the proof of Theorem \ref{main} is
almost identical to \cite[pp. 341-346]{Go-Hu1}, and so we shall
only sketch the remaining details.

Using Newton's identities and the polynomials above, we obtain the
same $\sigma_i's$ as in \cite[Eq. (18)]{Go-Hu1}, which shows that
$\sigma_i\in \mathbb{Z}[G]$ for $1\leq i\leq 5.$ Hence, from
\begin{equation} \label{Newt}
\phi^{\mu}|U_5 = \sigma_1\phi^{\mu-1}|U_5 -
\sigma_2\phi^{\mu-2}|U_5 + \sigma_3\phi^{\mu-3}|U_5 -
\sigma_4\phi^{\mu-4}|U_5 + \sigma_5\phi^{\mu-5}|U_5,
\end{equation}
for all $\mu\in Z,$ we obtain $\phi^{\mu}|U_5=\sum_{\nu=-\i}^\i
c_{\mu\nu}G^{\nu}$ with integral coefficients $c_{\mu\nu}$. And
arguing as in \cite[Lemma 6, Lemma 7]{Go-Hu1} gives
$$\pi(c_{\mu\nu})\geq [\frac{5\nu-3\mu-1}{6}]$$
from which we deduce Lemma \ref{a} by using
$G^i|U_5=(\phi^{2i}|U_5)G^{-i}$.
\end{proof}

\begin{lem}\label{b} For all $i$, $j\in
\mathbb{Z}$, we have $b_{ij}\in \mathbb{Z}$ and
$$\pi(b_{ij})\geq [\frac{5j-i-1}{6}]\quad \mathrm{and}
\quad \pi(b_{ij})\geq 1 \quad\mathrm{if}\quad i\equiv 1 \pmod
5.$$\end{lem} \begin{proof} We note that $F\phi^{\mu}|U_5$
satisfies the same Newton recurrence \eqref{Newt} as
$\phi^{\mu}|U_5$, and
$$F|U_5=5G,$$
$$F\phi^{-1}|U_5=-1,$$
$$F\phi^{-2}|U_5=q^{-1}+2+O(q)=G^{-1},$$
$$F\phi^{-3}|U_5=-5+O(q)=-5,$$
$$F\phi^{-4}|U_5=q^{-2}+9q^{-1}-9+O(q)=G^{-2}+5G^{-1}-25.$$
It hence follows that for all $\mu\in\mathbb{Z}$,
$$F\phi^{\mu}|U_5=\sum_{\nu}d_{\mu\nu}G^\nu,$$
where $d_{\mu\nu}$ are integers.  Arguing as in \cite[Lemma 8,
Lemma 9]{Go-Hu1}, we have that for all $\mu$, $\nu\in \mathbb{Z}$,
$$\pi(d_{\mu\nu})\geq [\frac{5\nu-3\mu-1}{6}]\quad \mathrm{and}
\quad \pi(d_{\mu\nu})\geq 1 \quad\mathrm{if}\quad \mu\equiv 2
\pmod 5.$$  Now, this lower bound for $\pi(d_{\mu\nu})$ along with
the fact $FG^i|U_5=(F\phi^{2i}|U_5)G^{-i}$ gives Lemma \ref{b}.
\end{proof}
Using the fact that $F|U_5 = 5G$ together with the lower bounds
for $\pi(a_{ij})$ and $\pi(b_{ij})$ in the lemmas above as in the
proof of \cite[Theorem 10]{Go-Hu1}, we obtain
\begin{lem}For all $\a\geq0$ and $j\geq 1$, we have
$$\pi(l_j(2\a+1))\geq\a+1+[\frac{j-1}{2}],$$
$$\pi(l_j(2\a+2))\geq\a+1+[\frac{j}{2}].$$\end{lem}
This implies that $l_j(2\a+1)\equiv 0 \pmod {5^{\a+1}}$, which is
Theorem \ref{main}.  This then completes the proof of Theorem
\ref{family}. \qed

For Theorem \ref{ramatype}, recall that we have shown that $F|U_5
= 5G$.  Hence we have
\begin{eqnarray*}
\frac{5\eta^3(\t)\eta(5\t)\eta^2(10\t)}{\eta^2(2\t)} &=&
\frac{\eta^3(\t)\eta(5\t)\eta^2(10\t)}{\eta^2(2\t)} \big | U_5 \\
&=& \left( (q^5;q^5)_{\infty}(q^{10};q^{10})_{\infty}^2
\sum_{n=0}^{\infty} \left(M_e(n)-M_o(n)\right)q^{n+1}\right) \big
| U_5 \\
&=& (q)_{\infty}(q^{2};q^{2})_{\infty}^2 \sum_{n=0}^{\infty}\left(
M_e(5n+4) - M_o(5n+4)\right)q^{n+1}.
\end{eqnarray*}
Multiplying the first and last terms in the above string of
equations by
$$
\frac{1}{q\prod_{n=1}^{\infty}(1-q^{n})(1-q^{2n})^2}
$$
gives Theorem \ref{ramatype}. \qed

\section{Proof of Theorem \ref{asy-cr}}\label{pr}

This proof is very similar to Kane's treatment of the circle
method in \cite{Kane} where he proved a conjecture by Andrews and
Lewis \cite{AL} on cranks of partitions modulo $3$.  By Cauchy's
integral formula, for a circle $C$ centered on the origin and
inside the unit circle, we have
\begin{eqnarray}\label{eq1}
M_e(n)-M_o(n)&=&\frac{1}{2\pi i}\int_C g(q)q^{-n-1}dq\cr
&=&\int_0^1 g\lt(\exp(-\rho+2\pi i\phi)\rt)e^{n\rho-2\pi i n
\phi}d\phi\cr &=& \mathop{\mathop{\sum_{0<k<N}}_{0<h<k}}_
{(h,k)=1}\int_{-\theta'_{h,k}}^{\theta''_{h,k}}g\lt(\exp\lt\{\frac{2\pi
i h}{k}-(\rho-2\pi i \phi)\rt\}\rt)e^{n(\rho-2\pi i
\phi)-\frac{2\pi i n h}{k}}d\phi,\end{eqnarray} where the last
equality follows from Andrews' dissection of the circle of
integration in \cite[Ch.~5]{Andrews} and $\rho$ is a positive real
number. As each of $\theta'_{h,k}$ and $\theta''_{h,k}$ is the
mediant of the Farey number $h/k$ and the adjacent Farey numbers,
it satisfies $1/(2kN)\leq\theta\leq 1/(kN)$.  If we substitute
$y=\rho-2\pi i \phi$ in (\ref{eq1}), we obtain
\begin{equation}\label{eq2}M_e(n)-M_o(n)=-\frac{1}{2\pi
i}\mathop{\mathop{\sum_{0<k<N}}_{0<h<k}}_ {(h,k)=1}\int_{\rho+2\pi
i\theta'_{h,k}}^{\rho-2\pi i\theta''_{h,k}}
g\lt(\exp\lt\{\frac{2\pi i h}{k}-y\rt\}\rt)e^{ny-\frac{2\pi i n
h}{k}}dy .\end{equation} Set
\begin{equation}\label{I}I:=-\frac{1}{2\pi i} \int_{\rho+2\pi
i\theta'_{h,k}}^{\rho-2\pi i\theta''_{h,k}}
g\lt(\exp\lt\{\frac{2\pi i h}{k}-y\rt\}\rt)e^{ny-\frac{2\pi i n
h}{k}}dy.\end{equation} Then we can write
\begin{equation}\label{sum}M_e(n)-M_o(n)=\mathop{\mathop{\sum_{2\nmid k,\
k<N}}_{0<h<k}}_ {(h,k)=1}I+\mathop{\mathop{\sum_{2|k,\
k<N}}_{0<h<k}}_ {(h,k)=1}I:=\Sigma_1+\Sigma_2.\end{equation}
Recall from \cite[Theorem 5.1]{Apos} that if $F(q)=1/(q;q)_\i$,
then
\begin{equation}\label{fe}
F\lt(\exp\lt(\frac{2\pi i h}{k}-\frac{2\pi z}{k^2}\rt)\rt)=e^{\pi
i
s(h,k)}\lt(\frac{z}{k}\rt)^{1/2}\exp\lt(\frac{\pi}{12z}-\frac{\pi
z}{12k^2}\rt)F\lt(\exp\lt(\frac{2\pi i
H}{k}-\frac{2\pi}{z}\rt)\rt),\end{equation} where $\Re(z)>0$,
$k>0$, $(h,k)=1$, and $hH\equiv -1\pmod k.$ We apply (\ref{fe}) to
the integrand of $I$ and evaluate a main term and estimate an
error term. We first find a bound for $\Sigma_1$:
\begin{lem}\label{sig1} If $\Sigma_1$ is as defined in (\ref{sum}), then
\begin{equation*}|\Sigma_1|<53n^{1/4}.\end{equation*}\end{lem}
\begin{proof}Since
\begin{equation}\label{g2}g(q)=(q;q)_\i(q;q^2)_\i^2
=\frac{(q;q)_\i^3}{(q^2;q^2)_\i^2}=\frac{F(q^2)^2}{F(q)^3},\end{equation}
if $k>0$ is not divisible by $2$, $h$, $H$ are integers defined as
in (\ref{fe}), $h'\equiv 2h \pmod k$, and $h'H'\equiv -1 \pmod k$,
then
$$g\lt(\exp\lt\{\frac{2\pi i h}{k}-y\rt\}\rt)
=2e^{\pi i
(2s(h',k)-3s(h,k))}\lt(\frac{2\pi}{yk}\rt)^{1/2}\exp\lt(-\frac{\pi^2}{3k^2y}-\frac{y}{24}\rt)
\frac{F^2\lt(\exp\lt\{\frac{2\pi i
H'}{k}-\frac{4\pi^2}{2k^2y}\rt\}\rt)}{F^3\lt(\exp\lt\{\frac{2\pi i
H}{k}-\frac{4\pi^2}{k^2y}\rt\}\rt)}.$$ Note that
$$\lt|\frac{F^2\lt(\exp\lt\{\frac{2\pi i
H'}{k}-\frac{4\pi^2}{2k^2y}\rt\}\rt)}{F^3\lt(\exp\lt\{\frac{2\pi i
H}{k}-\frac{4\pi^2}{k^2y}\rt\}\rt)}\rt|\leq
F^5\lt(\exp\lt\{-\frac{2\pi^2}{k^2}\Re(\frac{1}{y})\rt\}\rt)\leq
F^5\lt(\exp\lt\{-\frac{2\pi^2 a}{a^2+4\pi^2}\rt\}\rt),$$ because
\begin{equation}\label{rey}\Re(\frac{1}{y})=\frac{\rho}{\rho^2+4\pi^2\phi^2} \geq
\frac{k^2N^2\rho}{k^2N^2\rho^2+4\pi^2}\geq
\frac{k^2a}{a^2+4\pi^2}\quad \mathrm{if}\ a=N^2\rho.\end{equation}
We will choose $a$ in the interval $6.25\leq a\leq 6.5$ so that
$$F\lt(\exp\lt\{-\frac{2\pi^2 a}{a^2+4\pi^2}\rt\}\rt)\leq 1.35.$$
Also note that
$$\exp\lt(-\frac{\pi^2}{3k^2}\Re(\frac{1}{y})\rt)|y|^{-1/2}\leq\Re(y)^{-1/2}=\rho^{-1/2}$$
and the length of the integral is at most $4\pi/(kN)$.  Then,
\begin{eqnarray}\label{s1I}|\Sigma_1|&\leq& \mathop{\mathop{\sum_{2\nmid k,\
k<N}}_{0<h<k}}_ {(h,k)=1}\frac{1}{2\pi}\int_{\rho+2\pi
i\theta'_{h,k}}^{\rho-2\pi
i\theta''_{h,k}}\lt(\frac{8\pi}{k}\rt)^{1/2}
\lt|e^{(n-\frac{1}{24})y}\rt|\rho^{-1/2}
F^5\lt(\exp\lt\{-\frac{2\pi^2 a}{a^2+4\pi^2}\rt\}\rt)|dy| \cr
&\leq &4.4\mathop{\mathop{\sum_{2\nmid k,\ k<N}}_{0<h<k}}_
{(h,k)=1}\frac{2\pi}{k N}\lt(\frac{8\pi}{k}\rt)^{1/2}
e^{n\rho}\rho^{-1/2} \leq \ 45 \sum_{2\nmid k,\
k<N}\frac{e^{n\rho}\rho^{-1/2} }{k^{1/2}N}\cr &\leq &45
e^{n\rho}\rho^{-1/2}N^{-1/2}\leq \ 29
e^{n\rho}\rho^{-1/4}.\nonumber\end{eqnarray} Setting
$\displaystyle{\rho=\frac{1}{4n}}$, (that is, when $25n\leq
N^2\leq 26n$) yields Lemma \ref{sig1}.
\end{proof}

In order to compute $\Sigma_2$, consider $$\frac{1}{2\pi i}
\int_{\rho+2\pi i\theta'_{h,k}}^{\rho-2\pi i\theta''_{h,k}}
g\lt(\exp\lt\{\frac{2\pi i h}{2k}-y\rt\}\rt)e^{ny}dy.$$ Then by
the functional equation of the Dedekind-eta function in
(\ref{fe}), it is equal to \begin{equation}\label{I2}\frac{1}{2\pi
i} \int_{\rho+2\pi i\theta'_{h,k}}^{\rho-2\pi i\theta''_{h,k}}
\omega_{h,k}\sqrt{\frac{\pi}{yk}}e^{(n-1/24)y}\exp\lt(\frac{\pi^2}{24k^2y}\rt)g\lt(\exp\lt\{\frac{2\pi
i H}{2k}-\frac{4\pi^2}{4k^2y}\rt\}\rt)dy,\end{equation} where
$\omega_{h,k}= e^{\pi i (2s(h,k)-3s(h,2k))}$.  Now we approximate
$g$ by $1$. Let $\Psi(x)=g(x)-1$. Then
\begin{eqnarray}\label{s2}
\Sigma_2&=&\mathop{\mathop{\sum_{2k<N}}_{0<h<2k}}_
{(h,2k)=1}\frac{-1}{2\pi i}\exp\lt(-\frac{2\pi i n
h}{2k}\rt)\omega_{h,k}\sqrt{\frac{\pi}{k}} \int_{\rho+2\pi
i\theta'_{h,k}}^{\rho-2\pi i\theta''_{h,k}}
\sqrt{\frac{1}{y}}e^{(n-1/24)y}\exp\lt(\frac{\pi^2}{24k^2y}\rt)\cr
&&\qquad\qquad \times \lt\{\Psi\lt(\exp\lt\{\frac{2\pi i
H}{2k}-\frac{4\pi^2}{4k^2y}\rt\}\rt)dy+1dy\rt\}:=\Sigma_3+\Sigma_4.\end{eqnarray}
\begin{lem}\label{sig3} If $\Sigma_3$ is as defined in (\ref{s2}), then
\begin{equation*}|\Sigma_3|<22n^{1/4}.\end{equation*}\end{lem}
\begin{proof}
Since $|M_e(n)-M_o(n)|\leq p(n)$, $|g(x)-1|\leq|F(x)-1|,$ and
$$\lt|\frac{\Psi(x)}{x^{1/24}}\rt|\leq \frac{F(|x|)-1}{|x^{1/24}|}.$$ Thus
$$\frac{\lt|\Psi\lt(\exp\lt\{\frac{2\pi
i
H}{2k}-\frac{4\pi^2}{4k^2y}\rt\}\rt)\rt|}{\exp\lt(-\Re(1/y)\frac{\pi^2}{24k^2}\rt)}
\leq
\frac{F\lt(\exp\lt(-\Re(1/y)\frac{\pi^2}{k^2}\rt)\rt)-1}{{\exp\lt(-\Re(1/y)\frac{\pi^2}{24k^2}\rt)}}<
1.844,$$ by (\ref{rey}) and the choice of $a$. Hence
\begin{eqnarray*}\label{s3}|\Sigma_3|&<&1.844\mathop{\mathop{\sum_{2k<N}}_{0<h<2k}}_
{(h,2k)=1}\sqrt{\frac{1}{4\pi k}}\int_{\rho+2\pi
i\theta'_{h,k}}^{\rho-2\pi
i\theta''_{h,k}}\frac{1}{\Re(y)^{1/2}}\exp\lt(\Re(y)(n-1/24)\rt)|dy|\cr
&\leq&1.844\mathop{\mathop{\sum_{2k<N}}_{0<h<2k}}_
{(h,2k)=1}\sqrt{\frac{1}{4\pi k}}\frac{4\pi}{k
N}\rho^{-1/2}e^{\rho(n-1/24)}\leq
1.844\sum_{2k<N}\frac{4\sqrt{\pi}}{k^{1/2}N}\rho^{-1/2}e^{n\rho}\cr
&\leq&7.376\sqrt{2\pi}\rho^{-1/2}e^{n\rho}N^{-1/2}\leq
18.489\rho^{-1/4}e^{n\rho}a^{-1/4}.
\end{eqnarray*}
Since $\displaystyle{\rho=\frac{1}{4n}}$ and $6.25\leq a\leq6.5$,
we have Lemma \ref{sig3}.\end{proof}

Now, it remains to approximate $\Sigma_4$:
\begin{lem}\label{sig4} If $B_k(n)$ is as defined in Theorem \ref{asy-cr}, then
$$\Sigma_4=\frac{1}{\sqrt{n-1/24}}\sum_{0<k<5\sqrt{n}/2}
\frac{B_k(n)}{\sqrt{k}}\cosh\lt(\frac{\pi}{k}\sqrt{\frac{n-1/24}{6}}\rt)+E',$$
where $|E'|\leq 119n^{1/4}.$\end{lem}
\begin{proof}Consider
$$-\frac{1}{2\pi i}\int_{\rho+2\pi
i\theta'_{h,k}}^{\rho-2\pi i\theta''_{h,k}}
\sqrt{\frac{1}{y}}e^{(n-1/24)y}\exp\lt(\frac{\pi^2}{24k^2y}\rt)dy.$$
This is
$$\frac{1}{2\pi
i}\lt(\int_{-\i}^{(0+)}-\int_{-\i-2\pi i\theta''_{h,k}}^{\rho-2\pi
i\theta''_{h,k}}+\int_{-\i+2\pi i\theta'_{h,k}}^{\rho+2\pi
i\theta'_{h,k}}
\rt)\sqrt{\frac{1}{y}}e^{(n-1/24)y}\exp\lt(\frac{\pi^2}{24k^2y}\rt)dy:=J_0+J_1+J_2,$$
where $\int_{-\i}^{(0+)}$ denotes integration over the contour
leading from one branch of $-\i$ around $0$ to the other branch.
To compute an error contributed by $J_1$ and $J_2$, note that on
the lines $y=x+2\pi i \theta$, ($-\infty<x\leq \rho$,
$\theta=\pm\theta_{h,k}$ ) we have
$$\Re(\frac{\pi^2}{24k^2y})=\frac{x
\pi^2}{24k^2(x^2+4\pi^2\theta^2)}\leq
\frac{\rho}{96k^2\theta^2}\leq\frac{a}{24}\leq 0.271$$ and
$$|y|^{-1/2}=\lt(\frac{1}{\rho^2+4\pi^2\theta^2}\rt)^{1/4}\leq(\frac{1}{\theta})^{1/2}\leq\sqrt{2kN}.$$
Hence
$$|J_1+J_2|\leq \frac{2 e^{0.271}}{2\pi}\sqrt{2kN}\frac{e^{\rho(n-1/24)}}{n-1/24}
\leq \frac{e^{0.271}}{\pi}\sqrt{2kN}\frac{3e^{1/4}}{2n}\leq
1.137\sqrt{kN}/n.$$ Thus the total error made by $J_1$ and $J_2$
is at most
$$\mathop{\mathop{\sum_{2k<N}}_{0<h<2k}}_
{(h,2k)=1}1.137\frac{\sqrt{kN}}{n}\sqrt{\frac{\pi}{k}}\leq
1.137\frac{\sqrt{\pi N}}{n}\frac{N(N+1)}{2}\leq 2.016
N^{5/2}n^{-1}\leq 119n^{1/4}.$$ This leaves us to show that the
main term of $\Sigma_4$ is from $J_0$ which is equal to
$$
J_0=\frac{1}{2\pi
i}\int_{-\i}^{(0+)}\sqrt{\frac{1}{y}}e^{(n-1/24)y}\sum_{s=0}^\i
\lt(\frac{\pi^2}{24k^2}\rt)^s\frac{y^{-s}}{s!}dy.$$ Then by the
dominated convergence theorem and making a change of variable
$z=(n-1/24)y,$ we find that
\begin{eqnarray*}J_0&=&\sum_{s=0}^\i
\frac{1}{2\pi
i}\lt(\frac{\pi^2}{24k^2}\rt)^s\frac{1}{s!}\int_{-\i}^{(0+)}y^{-s-1/2}e^{(n-1/24)y}dy\cr
&=&(n-\frac{1}{24})^{-1/2}\sum_{s=0}^\i \frac{1}{2\pi
i}\frac{1}{s!}\lt(\frac{\pi^2(n-1/24)}{24k^2}\rt)^s\int_{-\i}^{(0+)}z^{-s-1/2}e^{z}dz.\nonumber\end{eqnarray*}
By Hankel's loop integral formula, we have
\begin{eqnarray*}J_0&=&(n-\frac{1}{24})^{-1/2}\sum_{s=0}^\i\frac{\lt(\frac{\pi^2(n-1/24)}{24k^2}\rt)^s}{s!\Gamma(s+1/2)}\cr
&=&\pi^{-1/2}(n-\frac{1}{24})^{-1/2}\sum_{s=0}^\i\frac{4^s\lt(\frac{\pi^2(n-1/24)}{24k^2}\rt)^s}{(2s)!}\cr
&=&\pi^{-1/2}(n-\frac{1}{24})^{-1/2}\cosh\lt(2\sqrt{\frac{\pi^2(n-1/24)}{24k^2}}\rt).\end{eqnarray*}
Therefore, we obtain the main term of $\Sigma_4$ as
$$\mathop{\mathop{\sum_{2k<N}}_{0<h<2k}}_
{(h,2k)=1}\omega_{h,k}\exp\lt(-\frac{2\pi i n
h}{2k}\rt)\sqrt{\frac{\pi}{k}}\pi^{-1/2}(n-\frac{1}{24})^{-1/2}\cosh\lt(2\sqrt{\frac{\pi^2(n-1/24)}{24k^2}}\rt).$$
Recalling $25n\leq N^2\leq 26n$ completes the proof of Lemma
\ref{sig4}.
\end{proof}
Theorem \ref{asy-cr} now follows from the three lemmas above.

\section{Weighted identities} In this section we prove Theorem
\ref{weighted} and discuss another weighted identity of similar
type.  We begin with a $q$-series of expansion of the crank
generating function \cite[Theorem 2.1]{BCCL}:
\begin{equation}\label{KW}
\frac{(q;q)_\i}{(xq;q)_\i(q/x;q)_\i}=
\frac{1}{(q;q)_{\infty}}\sum_{n=-\i}^\i\frac{(-1)^nq^{n(n+1)/2}(1-x)}{1-xq^n}.
\end{equation}
Setting $x=-1$, this may be rewritten in the following way:
\begin{equation} \label{KWnew}
\sum_{n \geq 0} (M_e(n) - M_o(n)) q^n = \frac{1}{(q;q)_{\infty}} +
4 \sum_{n \geq 1}
\frac{(-1)^nq^{n(n+1)/2}}{(q;q)_{n-1}(1-q^{2n})(q^{n+1};q)_{\infty}}.
\end{equation}

We shall interpret the right hand side of the identity as the
weighted count of partitions in Theorem \ref{weighted}.  How does
this go? The term $1/(q)_{\infty}$ initializes the weight of each
partition $\lambda$ to 1.  If there are no ones in $\lambda$, then
it is not counted at all by the sum on the right hand side, and so
the weight just remains 1.

Otherwise, we look at the ``initial run" of $\lambda$.  A
partition $\lambda$ will be counted by the sum on the right hand
side for each $j$ in the initial run that occurs an odd number of
times. For each such $j$, we add $4(-1)^j$ to the weight. (Those
that occur an even number of times contribute nothing.)  Then,
summing over all partitions of n, each counted according to its
weight, gives $M_e(n) - M_o(n)$.  \qed

As an application, we give a combinatorial proof of the following
$q$-series identity:
\begin{corol}
\begin{equation} \label{combproofeq}
\frac{1}{(q;q)_{\infty}} + 4\sum_{n=1}^{\infty}
\frac{(-1)^nq^{n(n+1)/2}}{(q;q)_{n-1}(1-q^{2n})(q^{n+1};q)_{\infty}}
= \sum_{n =0}^{\infty}
\frac{(-1)^nq^{n(n+1)/2}(1-q^{n+1})}{(q;q)_n(1+q^{n+1})(q^{n+2};q)_{\infty}}.
\end{equation}
\end{corol}
Of course, this identity may also be established analytically. For
example, take $a=b=0$, $c=d=q$, and $z=-1$ in the following
identity of S.H. Chan \cite[Eq. (3.1)]{Chan}, valid for $|a/c| <
1$ and $|bq/d| < 1$,
$$
\frac{(az,b/z,q)_{\infty}}{(cz,d/z)_{\infty}} =
\frac{(a/c,bc)_{\infty}}{(cd)_{\infty}}\sum_{n=0}^{\infty}
\frac{(cq/a,cd)_n(a/c)^n}{(q,bc)_n(1-czq^n)} +
\frac{d(ad,b/d)_{\infty}}{(cd)_{\infty}}\sum_{n=0}^{\infty}
\frac{(dq/b,cd)_n(bq/d)^n}{(q,ad)_n(z-dq^n)},
$$
where $(z_1,z_2,\dots ,z_k)_n=\prod_{i=1}^{k}(z_i;q)_n.$ The
result shows that the right hand side of \eqref{combproofeq} is
indeed the generating function for $M_e(n) - M_o(n).$ This with
(\ref{KWnew}) implies Corollary 4.1.
\begin{proof}
We define a second weight, $\omega_1$ by
$$
\omega_1(\lambda) = (-1)^{\text{length of initial run}}
-2\sum_{j}(-1)^j(-1)^{\text{$\#$ occurrences of $j$}},
$$
where the sum is over those $j$ in the initial run of $\lambda$.

We shall argue by induction on the length of the initial run that
for any partition $\lambda$, we have $\omega(\lambda) =
\omega_1(\lambda)$.  First, if the initial run is empty, then
$\omega(\lambda) = 1 = \omega_1(\lambda)$. Now, suppose $\lambda$
has an initial run of length $n+1$.  Let $\lambda -
\overline{n+1}$ denote the partition $\lambda$ with all of the
parts of size $n+1$ removed.  Then the induction hypothesis says
that $\omega(\lambda - \overline{n+1}) = \omega_1(\lambda -
\overline{n+1})$.

Let us now determine the effect on the weights of adding back in
the parts of size $n+1$.  First, consider the effect on $\omega$.
If $n+1$ occurs an even number of times, the weight is unchanged.
If it occurs an odd number of times, then $4$ is added to the
weight if $n$ is odd and $-4$ is added if $n$ is even.

Now, consider the effect on $\omega_1$.  If $n$ is even, the
length of the initial run changes from even to odd, giving us a
$-2$.  If $n+1$ occurs an even number of times, then we get a $+2$
and the net change is $0$.  If $n+1$ occurs an odd number of
times, we get a $-2$ and the net change is $-4$.  This matches the
change to $\omega$ in these cases.  A similar argument in the case
where $n$ is odd shows that the changes to $\omega$ and $\omega_1$
are always the same.  Hence, we have $\omega(\lambda) =
\omega_1(\lambda)$ for all partitions $\lambda$.

To complete the proof, it suffices to argue as in the proof of
Theorem \ref{weighted} that the right hand side of
\eqref{combproofeq} is the generating function for partitions
$\lambda$ counted with weight $\omega_1(\lambda)$.  The details
are very similar to the case of \eqref{KWnew}, so we omit them.
\end{proof}

Before continuing, we note that similar weighted identities can be
found for the rank difference $f(q)$, for example by using the
Watson's equation \cite{Watson}
\begin{equation}\label{W}
f(q)=\frac{1}{(q;q)_\i}\lt(1+4\sum_{k=1}^\i(-1)^k
\frac{q^{k(3k+1)/2}}{1+q^k}\rt).
\end{equation}

\section{Proof of Theorem \ref{formula}}\label{cr-dp}

In this section we prove Theorem \ref{formula}.  We begin by
deducing a key generating function for $M_e(\mathcal{D},n) -
M_o(\mathcal{D},n)$. It is easily seen using the definition of the
crank that
$$\sum_{n=1}^\i(M_e(\mathcal{D},n)-M_o(\mathcal{D},n))q^n=\sum_{n=1}^\i\frac{(-1)^{n+1} q^{n(n+3)/2}}{(-q;q)_n}
+\sum_{n=1}^\i\frac{(-1)^n q^{n(n+1)/2}}{(q;q)_{n-1}}.$$  But this
generating function does not tell us much about
$M_e(\mathcal{D},n) - M_o(\mathcal{D},n)$. However, we can prove
the more useful:
\begin{thm} \label{informative}
\begin{equation} \label{informativeeq}
\sum_{n=1}^\i(M_e(\mathcal{D},n)-M_o(\mathcal{D},n))q^n=\frac{1}{1+q}\sum_{n=1}^\i
q^{n(3n+1)/2}(1-q^{2n+1})-q(q^2;q)_{\infty}.
\end{equation}
\end{thm}

\begin{proof}
We shall demonstrate that
\begin{equation} \label{part1}
\sum_{n=1}^\i\frac{(-1)^{n+1} q^{n(n+3)/2}}{(-q;q)_n} =
\frac{1}{1+q}\sum_{n=1}^\i q^{n(3n+1)/2}(1-q^{2n+1})
\end{equation}
and
\begin{equation} \label{part2}
\sum_{n=1}^\i\frac{(-1)^n q^{n(n+1)/2}}{(q;q)_{n-1}} = -
q(q^2;q)_{\infty}.
\end{equation}
Shifting the summation variable by $1$ on both sides of
\eqref{part1}, we obtain the equivalent identity
\begin{equation}
\sum_{n = 0}^{\infty} \frac{(-1)^nq^{n(n+5)/2}}{(-q^2;q)_n} =
\sum_{n = 0}^{\infty} q^{(3n^2+7n)/2}(1-q^{2n+3}).
\end{equation}
But this is precisely the specialization $(a,b,c,d,e) \to
(q^3,-q^2,\infty,\infty,q)$ of the following limiting case of the
Watson-Whipple transformation:
\begin{equation} \label{Watson-Whipple}
\sum_{n=0}^{\infty}
\frac{(aq/bc,d,e)_n(\frac{aq}{de})^n}{(q,aq/b,aq/c)_n} =
\frac{(aq/d,aq/e)_{\infty}}{(aq,aq/de)_{\infty}}
\sum_{n=0}^{\infty}
\frac{(a,\sqrt{a}q,-\sqrt{a}q,b,c,d,e)_n(aq)^{2n}(-1)^nq^{n(n-1)/2}}
{(q,\sqrt{a},-\sqrt{a},aq/b,aq/c,aq/d,aq/e)_n(bcde)^n}.
\end{equation}
For \eqref{part2}, shifting the summation variable as before
leaves us with the task of proving that
\begin{equation}
\sum_{n=0}^{\infty} \frac{(-1)^nq^{n(n+3)/2}}{(q;q)_n} =
(q^2;q)_{\infty}.
\end{equation}
This follows immediately from the case $z=-q$ of
\begin{equation}
\sum_{n = 0}^{\infty} \frac{z^nq^{n(n+1)/2}}{(q;q)_n} =
(-zq;q)_{\infty}.
\end{equation}
This completes the proof Theorem \ref{informative}.
\end{proof}

We may now use this theorem to deduce the formula for
$M_e(\mathcal{D},n) - M_o(\mathcal{D},n)$.  First we record a
formula for each of the two series in \eqref{informativeeq}.
\begin{lem}
If $a(n)$ is defined by
\begin{equation} \label{aofnseries}
\sum_{n=1}^{\infty} a(n)q^n =
\frac{1}{1+q}\sum_{n=1}^\i q^{n(3n+1)/2}(1-q^{2n+1}),
\end{equation}
then
\begin{equation} \label{aofnformula}
a(n) =
\begin{cases}
-(-1)^{n - \lfloor n \rfloor_p}, & \text{if $R(\lfloor n
\rfloor_p)$ is even and positive}, \\
0, & \text{if $R(\lfloor n \rfloor_p)$ is odd and negative}, \\
(-1)^{n - \lfloor n \rfloor_p}, & \text{if $R(\lfloor n
\rfloor_p)$ is odd and positive}, \\
-2(-1)^{n - \lfloor n \rfloor_p}, & \text{if $R(\lfloor n
\rfloor_p)$ is even and negative}.
\end{cases}
\end{equation}
\end{lem}

\begin{proof}
The proof is elementary.  One simply expands $1/(1+q)$ as
$\sum(-1)^nq^n$ in \eqref{aofnseries}, multiplies the series
together, and verifies that the result is \eqref{aofnformula}.
\end{proof}

\begin{lem}
If $b(n)$ is defined by
\begin{equation} \label{bofnseries}
\sum_{n = 1}^{\infty} b(n)q^n = -q(q^2;q)_{\infty},
\end{equation}
then
\begin{equation} \label{bofnformula}
b(n) =
\begin{cases}
0, & \text{if $R(\lceil n \rceil_p)$ is positive}, \\
-1, & \text{if $R(\lceil n
\rceil_p)$ is odd and negative}, \\
1, & \text{if $R(\lceil n \rceil_p)$ is even and negative}.
\end{cases}
\end{equation}
\end{lem}

\begin{proof}
This is another elementary calculation, using the fact that
$$
(q^2;q)_{\infty} = \frac{1}{1-q}\sum_{n \in \mathbb{Z}}
(-1)^nq^{n(3n+1)/2}.
$$
\end{proof}

Theorem \ref{formula} now follows by combining the above two
lemmas and checking all of the different cases.    \qed

In closing this section, we might mention that it is also possible
to adapt Franklin's involution \cite[pp.10-11]{Andrews} to prove
Theorem \ref{formula}.

\section{Concluding Remarks}
We wish to end by offering two suggestions for future research.
First, there do not seem to be any simple congruences of the form
$M_e(pn+a) - M_o(pn+a) \equiv 0 \pmod{p}$ for $p$ prime except
when $p=5$.  Can the work of Ahlgren-Boylan \cite{AB} and
Kiming-Olsson \cite{KO} be adapted to prove that this is indeed
the case?  Second, can the ideas of \cite{LO} be applied to extend
the congruences modulo $5^a$ to congruences modulo $5^{a+1}$
within certain arithmetic progressions?

\end{document}